\documentclass[11pt]{amsart}
\usepackage{amsmath,amssymb}

\newtheorem{thm}{Theorem}
\newtheorem{lem}[thm]{Lemma}
\newtheorem{prop}[thm]{Proposition}

\newtheorem*{remark}{Remark}

\newcommand{\fq}{\mathbb{F}_q}

\newcommand{\F}{\mathcal{F}}

\numberwithin{equation}{section}

\newcommand{\Z}{{\mathbb Z}} 

\newcommand{\FF}{{\mathbb F}}

\newcommand{\E}{\mathbb{E}}
\newcommand{\prob}{\operatorname{Prob.}}

\newcommand{\sumpr}{\sideset{}{'}\sum}

\title[Number of points on hyperelliptic curves]{The fluctuations in
the number of points on a hyperelliptic curve over a finite field}

 \author{P\"ar Kurlberg}
\address{Department of Mathematics, Royal Institute of Technology,
 SE-100 44 Stockholm, Sweden}
 \email{kurlberg@math.kth.se}

\author{Ze\'ev Rudnick}
\address{Raymond and Beverly Sackler School of Mathematical Sciences,
Tel Aviv University, Tel Aviv 69978, Israel}
\email{rudnick@post.tau.ac.il}

\date{July 31, 2008}

 \thanks{P.K. was partially supported by grants from the G\"oran
 Gustafsson Foundation, the Knut and Alice Wallenberg foundation, the
 Royal Swedish Academy of Sciences, and the
 Swedish Research Council. Z.R. was supported by the Israel Science
Foundation (grant No. 925/06).}

\begin{document}
\begin{abstract}
The number of points on a hyperelliptic curve over a field of $q$
elements  may be expressed as $q+1+S$ where $S$ is a certain
character sum. We study fluctuations of $S$ as the curve varies
over a large family of hyperelliptic curves of genus $g$. For
fixed genus and growing $q$, Katz and Sarnak showed that
$S/\sqrt{q}$ is distributed as the trace of a random $2g\times 2g$
unitary symplectic matrix. When the finite field is fixed and the
genus grows, we find that the the limiting distribution of $S$ is
that of a sum of $q$ independent trinomial random variables taking
the values $\pm 1$ with probabilities $1/2(1+q^{-1})$ and the
value $0$ with probability $1/(q+1)$. When both the genus and the
finite field grow, we find that $S/\sqrt{q}$ has a standard
Gaussian distribution.
\end{abstract}

\maketitle

\section{Introduction}


Given  a finite field  $\fq$ of odd cardinality $q$ and
a square-free monic polynomial $F\in \FF_q[X]$  of degree $d\geq 3$,
we get  a smooth projective hyperelliptic curve $C_F$
with affine model
$$C_F:  Y^{2} = F(X)  $$
having genus $g=(d-2)/2$ when $d$ is even and
$g=(d-1)/2$ when $d$  is odd.
In this note we study the fluctuations in the number of $\FF_q$-points
on $C_F$ when  $F$ is drawn at random from the set 
of all  square-free monic polynomial $F\in \FF_q[X]$ of degree
$d$,  where the probability measure 
is obtained by picking the coefficients of $F$ uniformly in $\FF_q^d$ and
conditioning on $F$ being square free.
Correspondingly we get a probability measure on a family of
hyperelliptic  curves of genus $g\geq 1$ defined
over $\fq$. Our goal is to study these fluctuations
in the limit of either large genus or large $q$, or both.

%

The number of $\FF_q$-points on $C_F$ can be written
as\footnote{Giving the number of points of $C_{F}$ a spectral
  interpretation, it 
  is more natural to write the number of points as 
  $q-S'(F)+1$, where $S'(F)=-S(F)$ is the trace of the action induced
  by the Frobenius automorphism on a certain cohomology group.
  However, for our purposes, studying $S(F)$ will lead to slightly
  simpler notation, and, as we shall see, the distribution of $S(F)$ is
  symmetric, hence $S(F)$ and $S'(F)$ have the same distribution.}
$q+S(F)+1$ where $S(F)$ is the character sum 
$$
S(F) =\sum_{x \in \fq} \chi(F(x))
$$
and $\chi$ is the quadratic character of
$\FF_q^{\times}$ (with the convention that $\chi(0)=0$).  
Thus the problem is equivalent to
studying the fluctuations of $S(F)$ as $F$ varies over all square-free
polynomials in $\FF_q[X]$ of degree $d$, in the limit as either $d$ or
$q$ (or both) grow.

Our finding is that there are three distinct types of distribution
results according to the way the parameters $g$ and $q$ are allowed to
grow:

i) For $q$ fixed and the genus $g\to \infty$, we find that $S(F)$ is
distributed asymptotically as a sum of $q$ independent identically
distributed (i.i.d.) trinomial random variables $\{X_i\}_{i=1}^q$,
i.e., random variables taking values in $0,\pm 1$ with probabilities
$1/(q+1)$, $1/2(1+q^{-1})$ and $1/2(1+q^{-1})$ respectively.

ii) When  the genus $g$ is fixed and $q\to \infty$,
$S(F)/\sqrt{q}$ are
distributed as the trace of a random matrix in the group $USp(2g)$ of
$2g\times 2g$ unitary symplectic matrices. This is due to Katz and
Sarnak \cite{KatSar99}.

iii) When both $g\to \infty$ and $q\to \infty$ we find that
$S(F)/\sqrt{q}$ has a Gaussian value distribution with mean zero and
variance unity.

\bigskip

The case (iii) when both variables grow can be thought of as a
limiting case of either the two previous ones, when one of the two
parameters is held fixed. It is thus a good consistency check to see
that the limit distributions in both cases (i) and (ii) are a standard
Gaussian. Indeed, in the case when $q$ is fixed, (i) gives
that the limit distribution of $S(F)/\sqrt{q}$ is that of a
normalized sum $(X_1+\dots+X_q)/\sqrt{q}$ of $q$  i.i.d. random
variables; in turn the distribution of a normalized
sum of such i.id.'s converges, as $q\to\infty$, to a Gaussian
distribution with mean zero
and variance unity  by the Central Limit Theorem.
In case the genus $g$ is fixed, (ii) gives that the limit distribution of
$S(F)/\sqrt{q}$ is that of the traces of random matrices in
$USp(2g)$. The limit distribution of traces of a random matrix in
$USp(2g)$, as $g\to\infty$, is a standard Gaussian by a theorem
of Diaconis and Shahshahani \cite{DiaSha94}.
Of course this is not a proof of (iii), as it only addresses the
limiting form of the {\em limit distribution} in (i) and (ii), that is
either $\lim_{q\to\infty} (\lim_{g\to \infty})$ or
$\lim_{g\to\infty} (\lim_{q\to \infty})$
and not the joint limit  $\lim_{q\to\infty,g\to\infty}$.



\subsection{Some related work}

\bigskip
\noindent{\bf 1.}
In the unpublished manuscript \cite{Lar-gauss-preprint}, Larsen
studied moments for  a related family of hyperelliptic curves, namely
curves of the form $Y^{2} = \prod_{i=1}^{n}(X-a_{i})$, where $a_{1},
\ldots, a_{n}$ ranges over all $n$-tuples consisting of distinct
elements of $\fq$, and obtained Gaussian moments.

\bigskip
\noindent{\bf 2.}
Knizhnerman and Sokolinskii \cite{KniSok1979, KniSok1987}  computed
moments of the character sum $S(F)$ when $F$ ranges over all monic
non-square (rather than square-free) polynomials to investigate
extreme values taken by $S(F)$ (we thank Igor Shparlinski for this
reference).

\bigskip
\noindent{\bf 3.}
Bergstr\"{o}m \cite{Ber06} used methods closely
related to ours in order to obtain equivariant point counts for
families of hyperelliptic curves.  These point counts were then used
to determine cohomology groups of the moduli space of stable curves
of genus $2$ with $n$ marked points, for $n \leq 7$.

\bigskip
\noindent{\bf 4.} Finally, we refer to the recent preprint of 
Faifman and Rudnick \cite{FR08} which studies 
the statistics of the zeros of the zeta function of the curves $C_F$
over a fixed finite field in the limit of large genus.


\subsection{The main results}  Before giving a more quantitative
statement of our main results, we will need some notation.  Let
$V_{d} \subset \fq[X]$ be the set of monic polynomials of degree $d$,
and let $\F_{d} \subset V_{d}$ be the subset of {\em square free}
polynomials of degree $d$.
We will model $S(F)$ as a sum of $q$ independent identically
distributed (i.i.d.) trinomial random variables $\{X_i\}_{i=1}^q$,
where each $X_i$ takes values in $0,\pm 1$ with probabilities $1/(q+1)$,
$1/2(1+q^{-1})$ and $1/2(1+q^{-1})$ respectively.

For $q$ fixed and $d \to \infty$, we show that $S(F)$ behaves as
$\sum_{i=1}^q X_i$ in the
following sense:
\begin{thm}
\label{thm:distribution}
If $q$ is fixed and $d$ tends to infinity then the distribution of
$S(F)$, as $F$ ranges over all elements in $\F_d$,
is that of a sum of $q$ independent trinomial random variables.
More precisely, for $s \in \Z$ with $|s| \leq q$, we
have\footnote{Here, and in what follows, all  constants 
implied by the $O(\cdot)$-notation will be absolute.}
$$
\frac{|\{ F \in \F_d : S(F) = s  \}|}{|\F_{d}|}
= \prob
\left( \sum_{i = 1}^{q} X_{i} = s   \right)  \cdot (1+O(q^{(3q-d)/2})).
$$
%
\end{thm}

\begin{remark}
We may also let $q$ tend to infinity in
Theorem~\ref{thm:distribution}, provided that $d$ tends to
infinity in such a way that $d>3q$.
\end{remark}

By studying the moments we find that $S(F)/\sqrt{q}$ has
a Gaussian value distribution when both $d,q$  tend to infinity.
\begin{thm}
\label{thm:gaussian-moments}
  If $d,q$ both tend to infinity, then the moments of $S(F)/\sqrt{q}$ are
  asymptotically Gaussian with mean $0$ and variance $1$. 
In particular the limiting value distribution is a standard Gaussian.
\end{thm}

\section{Proof of Theorem~\ref{thm:distribution}}

The idea of the proof is to make the following heuristic precise:
Putting the uniform probability measure on $\F$, we may view $f \to
S(f)$ as a random variable on $\F$.  $S(f)$ can in turn be written as
$$
S(f) = \sum_{x \in \fq} X_{x}, 
$$
where for each $x \in \fq$,  $X_{x} = \chi(f(x))$ is also a random
variable on $\F$.  Then, 
as $d$ grows, the variables $\{X_{x}\}_{x \in \fq}$ become 
independent and the 
distribution of each individual $X_{x}$ is that of the earlier
mentioned trinomial random variable.

Thus, we will  study the following slightly more general problem: Given
a subset $S \subset \fq$ and a tuple $a = (a_x)_{x \in S}$, $a_x\in
\fq$, we wish to find the probability that for a randomly
selected $F \in \F$ we have $F(x) = a_{x}$ for all $x \in S$.

Before proceeding we need to introduce some additional notation.
For $F \in \fq[X]$, write $F = \prod_{i=1}^{n} F_{i}^{e_{i}}$ as a product
of irreducible polynomials, and let
$$
\mu(F) :=
\begin{cases}
0 & \text{if $e_{i}>1$ for some $i$,}\\
(-1)^{n} & \text{if $F$ is square free.}
\end{cases}
$$
Further, put
$$
|F| := q^{\deg(F)}
$$
and let
$$
\zeta_q(s) :=
\sum_{F \text{ monic}} |F|^{-s} =
\prod_{F \text{ irreducible}} (1-|F|^{-s})^{-1} = \frac{1}{1-q^{1-s}}
$$
be the (incomplete) zeta function of $\mathbb{A}^{1}/\fq$.

We will need to know the number of square free monic polynomials,
which can easily be deduced from the identity
$$ \zeta_q(s) = \zeta_q(2s) \sum_{d\geq 0} |\F_d| q^{-ds}, \qquad
\Re(s)>1.
$$
\begin{lem}
The number of square free monic polynomials of  degree $d$ equals
$$
|\F_d| =
\begin{cases}
q^{d}-q^{d-1} = q^d/\zeta(2) & \text{if $d \geq 2$,} \\
q^{d}  & \text{if $d = 0,1$.}
\end{cases}
$$
\end{lem}

We shall also need the following simple counting Lemma which
is at the heart of the independence result.
\begin{lem}
\label{l:surjectivity}
For $l \leq q$
let $x_1, x_2, \ldots, x_l \in \fq$ be distinct elements, and let
$a_1, a_2, \ldots, a_l \in \fq$.
If $d \geq l$, then
$$
|\{ F \in V_d : F(x_1) = a_1, \ldots, F(x_l) = a_l \}|
=
q^{d-l}.
$$
\end{lem}
%
%
%
%
%

\begin{proof}

Let $\tilde{V}_d = \{ g \in \fq[X] : \deg(g) \leq d-1\}$.  The map
$f(X) \to g(X) := f(X) - X^d$ then defines a bijection from $V_d$ to
$\tilde{V}_d$.
Since $f(x_i) = a_i$ for $1 \leq i \leq l$ is equivalent to $g(x_i) =
a_i - x_i^d$ for $1 \leq i \leq l$, we find that
\begin{multline}
\label{eq:counting-lemma}
|\{ f \in V_d : f(x_i) =
a_i \text{ for $1 \leq i \leq l$} \}|
\\ =
|\{ g \in \tilde{V}_d : g(x_i) = a_i-x_i^d \text{ for $1 \leq i \leq l$} \}|.
\end{multline}
Now,  the evaluation map $g \to (g(x_1), \ldots, g(x_l))$ is a
linear map from $\tilde{V}_d$ to $\fq^l$, and its kernel consists of
all $g \in \tilde{V}_d$ that are divisible by
$\prod_{i=1}^l(x-x_{i})$.  Hence the $\fq$-dimension of the kernel
equals $d-l$, and since $\dim_{\fq}(\tilde{V}_{d}) = d$ the  cokernel
has dimension $l$.  In particular, the
evaluation map is surjective, and both sides of
(\ref{eq:counting-lemma}) equal $q^{d-l}$ for all choices of $a_{1},
\ldots, a_{l}$.
\end{proof}

Next we determine the probability of a random polynomial in $\F_d$
taking a prescribed set of {\em nonzero} values on $l$ points.
\begin{lem}
\label{l:prob-of-nonzero-guys}
Let $d \geq 2$ and $l \leq q$ be a positive integers,
  let $x_1, x_2, \ldots, x_l \in \fq$ be distinct elements, and let
  $a_1, a_2, \ldots, a_l \in \fq$ be  {\em nonzero} elements.
Then
\begin{multline*}
\frac{|\{ F \in \F_d : F(x_1) = a_1, F(x_2)=a_2, \ldots, F(x_l) = a_l
  \}|}{|\F_d|}
\\=
\frac{q^{-l}}{(1-q^{-2})^{l}} \cdot
\left( 1+ O(q^{l-d/2}) \right).
\end{multline*}

\end{lem}
\begin{proof}
Using inclusion-exclusion, we find that
\begin{multline*}
|\{ F \in \F_d : F(x_i) = a_{i} \text{ for $1 \leq i \leq l$} \}|
=
\sum_{ \substack{ F \in V_{d} :
 F(x_i) = a_{i} \\ \text{ for $1 \leq i \leq l$}  }}
\mu(F)^{2}
\\
=
\sum_{D : \deg(D) \leq d/2}
\mu(D)
|\{ F \in V_{d-2\deg(D)} : D(x_i)^{2}F(x_i) = a_{i}
\text{ for $1 \leq i \leq l$}
\}|.
\end{multline*}

With $\sum^{'}$ denoting the sum over all polynomials such that $D(x)
\neq 0$ for all   $x \in \{x_1, x_2, \ldots, x_l\}$, we find, since
$a_{i} \neq 0$ for all
$i \leq l$, that the above equals
\begin{equation}
\label{eq:inc-excl}
\sumpr_{D : \deg(D) \leq d/2}
\mu(D)
|\{ F \in V_{d-2\deg(D)} : F(x_i) = a_{i} D(x_i)^{-2}
\text{ for $1 \leq i \leq l$}
  \}|.
\end{equation}
Now, as long as $\deg(F)=d-2\deg(D) \geq l$, by
Lemma~\ref{l:surjectivity},  we have
$$
|\{ F \in V_{d-2\deg(D)} :  F(x_i) = a_{i} D(x_i)^{-2}
\text{ for $1 \leq i \leq l$}
\}|
=
q^{d-2\deg(D)-l}
$$
hence (\ref{eq:inc-excl}) equals
\begin{equation}
\label{eq:main-and-error}
q^{d-l}
\sumpr_{D : \deg(D) < (d-l)/2}
\mu(D)
q^{-2\deg(D)}
+ \text{Error}
\end{equation}
where, since there can be at most one
polynomial $F$ of degree smaller
than $l$ that attains $l$ prescribed values (at $l$ distinct
points),
$$
\text{Error}
\ll
\sum_{D :  (d-l)/2 \leq \deg(D) \leq d/2 }
1
=
O(q^{d/2}).
$$

Our next goal is to evaluate the main term
$$
\sumpr_{D : 2\deg(D) < d-l}
\mu(D)
q^{-2\deg(D)}
=
\sumpr_{D}
\mu(D)
q^{-2\deg(D)}
+O(q^{(l-d)/2}).
$$
Noting that
\begin{multline*}
\sumpr_{D}
\mu(D)
|D|^{-2s}
=
\prod_{\substack{
F :  \text{ $F$ is irreducible,} \\ F(x_i) \neq 0 \text{ for $i \leq l$}}}
  (1-|F|^{-2s})
\\
=
(1-q^{-2s})^{q-l} \cdot
\prod_{\substack{\text{$F$ : $F$ is irreducible} \\ \deg(F)>1}}
(1-|F|^{-2s})
=
\frac{1}{\zeta(2s) (1-q^{-2s})^{l}}
\end{multline*}
we find that (\ref{eq:main-and-error}) equals
\begin{multline}
\label{counting-with-error-term}
q^{d-l}
\left(
\frac{1}{\zeta(2) (1-q^{-2})^{l}}
+ O(q^{(l-d)/2})
\right)
+ O(q^{d/2})
=
\frac{q^{d-l}}{\zeta(2) (1-q^{-2})^{l}}
+ O(q^{d/2}).
\end{multline}
Since $|\F_d| = \frac{q^d}{\zeta(2)}$ for $d \geq 2$,  the probability
that $F(x_i)=a_{i}$ for all $i \leq l$ equals
$$
\frac{q^{-l}}{(1-q^{-2})^{l}}
+ O(q^{-d/2}),
$$
concluding the proof.
\end{proof}

We now easily obtain the probability of $F$ attaining any set of
prescribed values.
\begin{prop}
\label{prop:probability-with-zero-guys}
Let $x_{1}, x_{2}, \ldots, x_{l}, x_{l+1}, x_{l+m} \in \fq$ be
distinct elements, let $a_{1}, a_{2}, \ldots a_{l} \in \fq^{\times}$,
and let $a_{l+1}=a_{l+2}= \ldots = a_{l+m} = 0$.
Then
\begin{multline}
\frac{1}{|\F_{d}|}
|\{ F \in \F_{d} : F(x_{i}) = a_i \text{ for $1 \leq i \leq m+l$} \}|
\\ =
\frac{(1-1/q)^{m} q^{-(m+l)} }{(1-q^{-2})^{m+l}} \cdot
\left(1+O(q^{(3m+2l-d)/2})\right).
\end{multline}
\end{prop}

\begin{proof}
Any  $F \in \F_d$ which vanishes at $Z = \{x_{l+1}, x_{l+2}, \ldots,
x_{l+m}\}$ can be written as
$$
F(x) = \prod_{i=l+1}^{l+m}(x-x_{i}) G(x)
$$
where $G \in \F_{d-m}$ is a square free polynomial that is
non-vanishing on $Z$.
Moreover, the condition that $F(x_i) = a_{i}$ for $1 \leq i \leq l$
can then be expressed as $G(x_i) = a_{i}
\prod_{j=l+1}^{l+m}(x_i-x_{j})^{-1}$ for $1 \leq i \leq l$, and $G(x_j)$
is arbitrary (but nonzero) for $l+1 \leq j \leq l+m$.  In other words,
there are $(q-1)^{m}$ possible values for $G$ restricted to $Z$ and by
the Lemma~\ref{l:prob-of-nonzero-guys} (in particular,
see~(\ref{counting-with-error-term})), the number of such polynomials
equals
\begin{multline*}
(q-1)^{m} \left(
\frac{q^{d-m-(m+l)}}{\zeta(2) (1-q^{-2})^{m+l}}
+O(q^{(d-m)/2})
\right)
\\ =
(1-1/q)^{m} \left(
\frac{q^{d-(l+m)}}{\zeta(2) (1-q^{-2})^{m+l}}
+O(q^{(d+m)/2})
\right).
\end{multline*}

Dividing by the number of square free
polynomials, we find that the probability of a random $F \in
  \F_d$ vanishing on $Z$, and taking prescribed values outside $Z$
equals
$$
\frac{(1-1/q)^{m} q^{-(m+l)}}{(1-q^{-2})^{m+l}} (1 + O(q^{(3m+2l-d)/2})).
$$
\end{proof}

To finish the proof of Theorem~\ref{thm:distribution}, we argue as
follows:
Let $x_{1}, x_{2}, \ldots, x_{q}$ be distinct elements of $\fq$,
let $\epsilon_i \in \{-1,0,1\}$ for $1 \leq i \leq q$, and
define $m=|\{ i \in \{1,2,\ldots,q\}: \epsilon_i = 0\}|$.
Taking $l=q-m$ in Proposition~\ref{prop:probability-with-zero-guys}
and noting that the number of nonzero squares, respectively
non-squares, in $\fq$ equals $(q-1)/2$, we find that
\begin{multline*}
\frac{
|\{ F \in \F_d : \chi(F(x_{i})) =
\epsilon_i \text{ for all $1 \leq i \leq q$}  \}|
}{|\F_d|}
\\=
\left( \frac{q-1}{2} \right)^{q-m} \cdot
\frac{(1-1/q)^{m} q^{-q}}{(1-q^{-2})^{q}} (1 + O(q^{(3q-d)/2}))
\\=
2^{-(q-m)}
\frac{(1-1/q)^{q} q^{-m}}{(1-q^{-2})^{q}} (1 + O(q^{(3q-d)/2}))
\\=
\frac{2^{-(q-m)} q^{-m}}{(1+q^{-1})^{q}} (1 + O(q^{(3q-d)/2})).
\end{multline*}
On the other hand, if $X_{1}, \ldots, X_{q}$ are i.i.d.
trinomial random variables as before, we have
\begin{multline*}
\prob( X_{i} = \epsilon_{i} \text{ for $1 \leq i \leq q$})
=
(q+1)^{-m} \cdot 2^{-(q-m)} (1+q^{-1})^{m-q}
\\ =
2^{-(q-m)}    (1+q^{-1})^{-q} q^{-m}.
\end{multline*}
Summing
over all possible choices of $\{\epsilon_i\}_{i=1}^{q}$ such that
$\sum_{i=1}^{q} \epsilon_i = s$, the proof is concluded.

\section{Proof of Theorem~\ref{thm:gaussian-moments}}
\label{sec:moments}

Let
$$
M_k(q,d) :=
\frac{1}{|\F_d|}
\sum_{F \in \F_d}  \left( \frac{S(F)}{\sqrt{q}} \right)^k
$$
be the $k$-th moment of $S(F)$ as $F$ ranges over the family of square
free polynomials of degree $d$ in $\fq[X]$.  As before, let $X_{1},
\ldots, X_{q}$ be independent trinomial random variables, taking
values $-1,0,1$ with probabilities $(\frac{q/2}{q+1},\frac{1}{q+1}
,\frac{q/2}{q+1})$.  Theorem~\ref{thm:gaussian-moments} is then an
immediate consequence of the following Proposition.
\begin{prop}
We have
$$
M_{k}(q,d)
=
\E \left(
\left(
\frac{1}{q^{1/2}}
\sum_{i=1}^q X_{i}
\right)^{k} \right)
+
O(q^{(3k-d)/2}).
$$
In particular, if
$q,d \to \infty$,
$M_k(q,d)$ agrees with Gaussian moments for all $k$.
\end{prop}

\begin{proof}  We have
  \begin{multline}
\label{eq:moment}
M_k(q,d) =
\frac{1}{|\F_d|}
\sum_{f \in \F_d}
\left( \frac{1}{q^{1/2}}
\sum_{x \in \fq}  \chi(f(x)) \right)^k
\\=
\frac{1}{q^{k/2}}
\sum_{x_{1}, x_{2}, \ldots, x_{k} \in \fq}
\sum_{f \in \F_d}
\chi(f(x_{1})f(x_{2}) \cdots f(x_{k}))
\\
=
\frac{1}{q^{k/2}}
\sum_{l=1}^k
c(k,l)
\sum_{((x_{1}, \ldots, x_{l}), (\epsilon_1, \ldots,\epsilon_{l})) \in P_{k,l}}
\frac{1}{|\F_d|}
\sum_{f \in \F_d}
\chi
\left( \prod_{i=1}^l f(x_{i})^{\epsilon_i} \right)
\end{multline}
where
$$
P_{k,l} = \{ ((x_{1}, \ldots, x_{l}), (\epsilon_1, \ldots,
\epsilon_{l}))
: x_{1}, \ldots, x_{l} \text{ all distinct and } \sum_{i=1}^l \epsilon_i = k
\}
$$
and $c(k,l)$ is a certain combinatorial factor, whose exact form is
unimportant.
Now, by Lemma~\ref{l:prob-of-nonzero-guys},
$$
\frac{1}{|\F_d|}
\sum_{f \in \F_d}
\chi
\left( \prod_{i=1}^l f(x_{i})^{\epsilon_i} \right)
= 0 + O( q^{l-d/2})
$$
unless {\em all} $\epsilon_i$ are even, in which case
$$
\frac{1}{|\F_d|}
\sum_{f \in \F_d}
\chi
\left( \prod_{i=1}^l f(x_{i})^{\epsilon_i} \right)
= \frac{1}{(1+1/q)^{l}} +  O( q^{l-d/2}).
$$
Noting that
$$
\sum_{l=1}^k
c(k,l)
\sum_{((x_{1}, \ldots, x_{l}), (\epsilon_1, \ldots,\epsilon_{l})) \in P_{k,l}} 1
= q^k
$$
we find that the contribution from the error terms is $\ll q^{-k/2}
q^{k} q^{k-d/2}$ and thus (\ref{eq:moment}) equals
$$
\frac{1}{q^{k/2}}
\sum_{l=1}^k
c(k,l)
\sum_{\substack{ ((x_{1}, \ldots, x_{l}), (\epsilon_1,
    \ldots,\epsilon_{l})) \in   P_{k,l} \\ \text{ all $\epsilon_i$ even}}}
\frac{1}{(1+1/q)^{l}}
+
O(  q^{(3k-d)/2} ).
$$

On the other hand, since $X_{1}, \ldots, X_{q}$ are independent
trinomial random variables, we find that the
expectation of the $k$-th moment of their normalized sum is
\begin{multline*}
\E \left( \left(
\frac{1}{q^{1/2}}
\sum_{j=1}^{q} X_{j}
\right)^{k}\right)
=
\frac{1}{q^{k/2}}
\sum_{i_{1}, i_{2}, \ldots, i_{k} \in \{0,1,\ldots,q-1\}}
\E (X_{i_{1}} \cdot X_{i_{2}} \cdots X_{i_{k}})
\\=
\frac{1}{q^{k/2}}
\sum_{l=1}^k
c(k,l)
\sum_{((x_{1}, \ldots, x_{l}), (\epsilon_1, \ldots,\epsilon_{l})) \in P_{k,l}}
 \E \left( \prod_{j=1}^l X_{i_{j}}^{\epsilon_i} \right).
\end{multline*}
As before we have $\E \left( \prod_{j=1}^l
  X_{i_{j}}^{\epsilon_i} \right) = 0$,
unless all $\epsilon_i$ are even, in which case
$$
\E \left( \prod_{j=1}^l  X_{i_{j}}^{\epsilon_i} \right)  =
\frac{1}{(1+1/q)^{l}}
$$
(note that $\E(X_{i_j}^\epsilon) = 0$ for $\epsilon$ odd, and
$\E(X_{i_j}^{\epsilon}) = 1/(1+q^{-1})$ for $\epsilon$ even),
concluding the proof of the first  assertion.

The final assertion follows since the moments of a sum of
bounded i.i.d. random variables converge to the Gaussian moments,
cf. \cite[section 30]{bil95}.
\end{proof}



\providecommand{\bysame}{\leavevmode\hbox to3em{\hrulefill}\thinspace}
\providecommand{\MR}{\relax\ifhmode\unskip\space\fi MR }
\providecommand{\MRhref}[2]{%
  \href{http://www.ams.org/mathscinet-getitem?mr=#1}{#2}
}
\providecommand{\href}[2]{#2}

\end{document}